\documentclass[12pt, a4paper]{article}
\usepackage[margin=3.5cm]{geometry}

\usepackage{amsmath,amssymb,amsthm}
\usepackage{algorithm} 
\usepackage[noend]{algpseudocode}

\usepackage[show]{ed}

\usepackage{tikz}

\newtheorem{lemma}{Lemma}[section]

\newtheorem{theorem}[lemma]{Theorem}

\newcommand{\beq}{\begin{equation}}
\newcommand{\eeq}{\end{equation}}

\makeatletter
\newcommand{\bigger}{\bBigg@{3}}
\newcommand{\vast}{\bBigg@{4}}
\newcommand{\Vast}{\bBigg@{5}}
\makeatother

\makeatletter
\@tfor\next:=abcdeghijklmnopqstuvwyxABCDEFGHIJKLMNOPQRSTUVWXYZ\do{%
  \def\command@factory#1{%
    \expandafter\newcommand\csname b#1\endcsname{\boldsymbol{#1}}
    \expandafter\newcommand\csname c#1\endcsname{\mathcal{#1}}
    \expandafter\newcommand\csname bb#1\endcsname{\mathbb{#1}}
  }
 \expandafter\command@factory\next}

\def\greekvectors#1{%
 \@for\next:=#1\do{%
    \def\X##1;{%
     \expandafter\newcommand\csname b##1\endcsname{\boldsymbol{\csname##1\endcsname}}
     }
   \expandafter\X\next;
  }
}
\greekvectors{alpha,beta,gamma,Gamma,delta,Delta,epsilon,varepsilon,zeta,theta,vartheta,Theta,iota,kappa,lambda,Lambda,mu,nu,xi,Xi,pi,Pi,rho,varrho,sigma,Sigma,tau,upsilon,Upsilon,phi,varphi,Phi,chi,psi,Psi,omega,Omega} 

\def\<{\langle}
\def\>{\rangle}
\newcommand\e{\varepsilon}
\def\wt{\widetilde}

\def\Chi{\raise .3ex
\hbox{\large $\chi$}}

\newcommand{\spn}{\mathop{\operator@font span}} 
\newcommand{\cspn}{\overline{\mathop{\operator@font span}}}
\newcommand{\dist}{\mathop{\operator@font dist}}
\newcommand{\argmin}{\mathop{\operator@font argmin}}
\newcommand{\argmax}{\mathop{\operator@font argmax}}
\newcommand{\range}{\mathop{\operator@font Ra}}
\newcommand{\image}{\mathop{\operator@font Im}}

\title{Coarse reduced model selection for nonlinear state estimation}
\author{J.A.~Nichols}
\date{\today}

\begin{document}
\maketitle
\begin{abstract}
State estimation is the task of approximately reconstructing a solution \(u\) of a parametric partial differential equation when the parameter vector \(y\) is unknown and the only information is \(m\) linear measurements of \(u\). In \cite{cohen_nonlinear_2020} the authors proposed a method to use a family of linear reduced spaces as a generalised nonlinear reduced model for state estimation. 
A computable {\em surrogate distance} is used to evaluate which linear estimate lies closest to a true solution of the PDE problem. 

In this paper we propose a strategy of coarse computation of the surrogate distance while maintaining a fine mesh reduced model, as the computational cost of the surrogate distance is large relative to the reduced modelling task. We demonstrate numerically that the error induced by the coarse distance is dominated by other approximation errors.
\end{abstract}

\section{Introduction}
%
Complex physical systems are often modelled by a parametric partial differential equation (PDE). We consider the general problem of the form
\begin{equation} \label{eq:pde}
A(y) u = f(y),
\end{equation}
where $A(y) : \cH \to \cH^\prime$ is an elliptic second order operator, $\cH$ is an appropriate Hilbert space. The problem is defined on a physical domain $D \subset \bbR^{2,3}$, and the parameter $y$ is within a parameter domain $\cY$, typically a subset of $\bbR^d$ where $d$ might be large. Each parameter $y\in \cY$ is assumed to result in a unique solution $u(x,y)$ but we mostly suppress the spatial dependence $x \in D$ and write $u(y)$. We denote all $\cH$-norms by $\| \cdot \| := \| \cdot \|_{\cH}$ and inner-products $\<\cdot,\cdot\> := \< \cdot, \cdot\>_{\cH}$, and use an explicit subscript if they are from a different space.

In practical modelling applications it is often computationally expensive to produce a high precision numerical solution to the PDE problem \eqref{eq:pde}. To our advantage however, the mapping $u: y \to \cH$ is typically smooth and compact, hence the set of solutions over all parameter values will be a smooth manifold, possibly of finite intrinsic dimension. We define the {\em solution manifold} as follows, assuming from here that it is compact,
\[
\cM := \{ u(y) : y \in \cY\} .
\]
%

The methods proposed in this paper extend on {\em reduced basis} approximations. A reduced basis is a linear subspace $V\subset \cH$ of moderate dimension $n := \dim(V) \le \dim(\cM) \le \dim(\cH)$. We use the {\em worst-case error} as a benchmark, defined as
\begin{equation} \label{eq:wc_err}
\e(V, \cM) = \max_{u \in \cM} \dist(u, V) = \max_{u \in \cM} \| u - P_V u \|,
\end{equation}
where $P_V$ is the orthogonal projection operator on to $V$.
There are a variety of methods to construct a reduced basis with desirable worst-case error performance, and here we concentrate on {\em greedy methods} that select points in $\cM$ which become the basis for $V$. We discuss these methods further in \S \ref{sec:lin_methods}.

A reduced basis can be used to accelerate the {\em forward problem}. One can numerically solve the PDE problem for a given parameter $y$ by using $V$ directly in the Galerkin method, making the numerical problem vastly smaller while retaining a high level of accuracy. 
A thorough treatment of the development of reduced basis approximations is given in \cite{hesthaven_certified_2016}.

In this paper we are concerned with {\em inverse problems}. In this setting it is assumed that there is some unknown true state $u$ (which could correspond to the state of some physical system), and we do not know the parameter vector $y$ that gives this solution. Instead, we make do with a handful of $m$ linear measurements $\ell_i(u)$. These measurements are used to make some kind of accurate reconstruction of $u$ ({\em state estimation}) or a guess of the true parameter $y$ ({\em parameter estimation}).

The {\em Parametrized Background Data Weak} (PBDW) approach introduced in \cite{maday_parameterized-background_2015} gives a straightforward procedure for finding an estimator $u^*$ of the true state $u$, using only the linear measurement information and a reduced basis $V$. One limitation however is that the estimator error $\| u^* - u\|$ is bounded from below by the {\em Kolmogorov $n$-width} given by
\begin{equation} \label{eq:kol_n}
d_n(\cM) = \inf_{\dim(V) = n} \sup_{u \in \cM} \dist(u, V),
\end{equation}
 where the infimum is taken over all $n$ dimensional linear spaces in $\cH$. The $n$-width is known to converge slowly for many parametric PDE problems.
 
We will review methods for constructing a family of {\em local linear} reduced models and a nonlinear estimator $u^*$ using a surrogate distance model selection procedure. We propose the use of coarser finite element meshes to perform this selection. This coarse selection strategy is motivated by the observation in \cite{cohen_nonlinear_2020} that the model selection is by far the most computationally costly component in the nonlinear estimation routine. Finally in \S\ref{sec:numerics} we examine numerical examples of surrogate distances over different mesh widths, and see that they make insubstantial impacts to the nonlinear estimator.

\subsection{Linear reduced models} \label{sec:lin_methods}

A reduced basis is a linear space of the form $V = \spn( u^1, \ldots, u^n ) \subset \cH$, where the $u^i = u(y^i) \in \cM$. The parameter values $y^i$ are typically chosen in some iterative greedy procedure to try and minimise $\e(V, \cM)$ at each step.

We can define a greedy procedure follows: given $V$ of dimension $n$ and  a finite subset $\wt \cY$ of $\cY$, to produce an $(n+1)$-dimensional reduced space we find the parameter $y^{n+1}\in \wt \cY$ which gives us the largest  $\dist(u( y^{n+1}), V) = \| u( y^{n+1}) - P_V u( y^{n+1}) \|$. We then augment the space $V$ with $u^{n+1} = u( y^{n+1})$. This simple strategy, in some cases \cite{binev_convergence_2011}, yields a reduced basis that is optimal with regards to the Kolmogorov $n$-width of $\cM$.
In this setting the quantity 
\begin{equation} \label{eq:err_est}
\e_{\mathrm{est}}(V, \cM) := \max_{y\in \wt\cY} \| u(y) - P_V u(y) \|
\end{equation} serves as a reasonable and calculable estimate of $\e(V, \cM)$, as shown in \cite{cohen_reduced_2020}.

In practice, we also allow $V$ to be an {\em affine} space, with an offset $\bar u$, such that
\[
V = \bar u \oplus \spn(u^1,\ldots,u^n) .
\]
Typically we take $\bar u$ to be an approximate barycenter of $\cM$. From here we will use the term {\em reduced space} or {\em model} rather than basis.

With $V$ chosen, we now consider the state estimation problem. We have $m$ pieces of linear data $\ell_i(u)$ for $i=1,\ldots, m$ of the unknown state $u$, and the $\ell_i \in \cH^\prime$. We assume we know the {\em form} of the functionals $\ell_i$ and hence the Riesz representers $\omega_i$ for which $\< \omega_i, u \> = \ell_i(u)$. These $\omega_i$ define a {\em measurement space} and the measurement vector of $u$:
\[
W := \spn(\omega_1, \ldots, \omega_m) \quad\text{and}\quad w = P_W u.
\]
Note here that we assume no noise in our measurements, but allowing for random noise is straightforward and has been considered in \cite{maday_parameterized-background_2015, binev_data_2017}.

The PBDW approach, developed in \cite{maday_parameterized-background_2015}, seeks a reconstruction candidate or estimator $u^*(w)$ that is close to $u$, but that agrees with the measurement data, that is $P_W u^*(w) = P_W u = w$. They define an estimator
\[
u^*(w) = \argmin \{ \dist(v, V) : P_W v = w \},
\]
which can be calculated through a set of normal equations of size $n\times m$, using the cross-Grammian matrix of the bases of $W$ and $V$.
Given only the measurement information $w$, the measurement space $W$, and the reduced space $V$, this estimator $u^*(w)$ is an optimal choice \cite{binev_data_2017}. 
The estimator lies in the subspace $u^*(w) \in V \oplus W$.

We require that $W^\perp \cap V = \emptyset$ for this reconstruction algorithm to be well posed (otherwise there are infinitely many candidates for $u^*$), which in turn requires that 
$n = \dim(V) < \dim(W) = m$. This dimensionality requirement is reflected in the error analysis. We define an inf-sup constant
\begin{equation} \label{eq:inf_sup}
\mu(V, W) := \max_{v \in V} \frac{ \| v \|} {\| P_W v \|}
\end{equation}
which is the inverse of the cosine of the angle between $V$ and $W$ and we have $\mu(V, W) \in [1, \infty]$. For $\mu$ to be finite, we require $W^\perp \cap V = \emptyset$.
The inf-sup constant plays the role of a stability constant for our linear estimator as we have the well known bound 
\[
E_{\mathrm{wc}} = \max_{u \in \cM} \| u - u^*(P_{W_m} u) \| \le \mu(V, W) \, \e(V, \cM) ,
\]
as deomstrated in \cite{binev_data_2017}.
From the definitions we have $\e(V, \cM) \ge d_{n+1}(\cM)$, hence this reconstruction error can at best be the $(n+1)$-width of $\cM$.

\section{Nonlinear reduced models} \label{sec:nonlinear}

A fundamental drawback of linear reduced models is the slow decay of the Kolmogorov $n$-width for a wide variety of PDE problems. To circumvent this limitation, a framework for non-linear reduced models and their use for state estimation was presented in \cite{cohen_nonlinear_2020}. The proposal involves determining a partition of the manifold
\[
\cM = \bigcup_{k=1}^K \cM_k
\]
and producing a family of affine reduced space approximations $V_k$ to each portion $\cM_k$. Each space has dimension $n_k = \dim(V_k)$, requiring $n_k < m$ for well-posedness. 

Given any target $\e>0$ and $\mu\ge 1$ it is possible with large enough $K$ to determine a partition and family of reduced spaces $V_k$  that satisfies
\begin{equation} \label{eq:mu_eps_adm}
\e_k := \e(V_k, \cM_k) \le \e \quad \text{and} \quad \mu_k := \mu(V_k, W) \le \mu \quad \text{for all } k=1,\ldots,K
\end{equation}
in which case we say the family $(V_k)_{k=1}^K$ is $(\mu,\e)$-admissible. 
A slightly looser criteria on the partition can also be satisfied: given some $\sigma>0$, we say the family $(V_k)_{k=1}^K$ is {\em $\sigma$-admissible} if $\mu_k \e_k \le \sigma$ for all $k=1,\ldots,K$.
The existence of both $(\mu,\e)$ and $\sigma$-admissible families follow from the compactness of $\cM$, and a full demonstration can be found in \cite{cohen_nonlinear_2020}. 

In practice one may construct a $\sigma$-admissible family in the following way. Say we are given a partition of the parameter space $(\cY_k)_{k=1}^{K-1}$ with $\cup_{k=1}^{K-1} \cY_k = \cY$, and the associated partition of the manifold $(\cM_k)_{k=1}^{K-1}$. We have a reduced space $V_k$, produced by a greedy algorithm on each $\cM_k$. With each $V_k$ we have an associated error estimate $\sigma_{\mathrm{est, k}} = \mu_k(V, W) \e_{\mathrm{est}, k}(V, \cM)$. We can pick the the largest $\sigma_{\mathrm{est, k}}$, with index $\tilde k$ say, and we split the cell $\cY_{\tilde k}$ in half for each parameter coordinate direction $i \in \{1,\ldots,d\}$, resulting in two reduced spaces $V_{ \tilde k, i}^+$ and $V_{ \tilde k, i}^-$ for each split direction. We take the split direction $i$ to be the one with the smallest maximum error $ \max( \sigma_{\tilde k, i}^+, \sigma_{\tilde k, i}^-)$, and we enrich the family $(V_k)_{k=1}^{K-1}$ with the two new reduced spaces, making sure to remove $V_{\tilde k}$ from the collection. More details of the splitting procedure can be found in \cite{cohen_nonlinear_2020}.

\subsection{Surrogate reduced model selection} \label{sec:selec}

Based on a measurement $w$, each affine reduced space has an associated reconstruction candidate that is found through the PBDW method,
\begin{equation}
u^*_k(w) := \argmin \{ \dist(v, V_k) : P_W v = w \} \text{ for } k=1,\ldots,K
\end{equation}
If we happen to know which that the true state $u$ originates from some $\cM_k$, then we would best use $u^*_k$ as our estimator. In this scenario we would have an error bound of $\| u - u^*_k(P_W u) \| \le \e_k \mu_k < \sigma$. This information about the true state $u$ is not available in practice, so we require some other method to determine which candidate $u^*_k$ to choose.

Consider a surrogate distance $\cS(v, \cM)$ from $v$ to $\cM$ that satisfies the uniform bound for $0 < r \le R $,
\begin{equation} \label{eq:surr_bounds}
r \dist (v, \cM) \le \cS(v, \cM) \le R \dist (v, \cM).
\end{equation}
If this surrogate distance is computable, then we can use it to find a {\em surrogate selected} estimator by choosing
\begin{equation} \label{eq:surr_sel}
k^* := \argmin \{ \cS(u^*_k, \cM) : k=1,\ldots,K \},\text{ and define } u^*(w) := u^*_{k^*}(w),
\end{equation}
noting that as there is a dependence on $w$ we can write $k^*(w)$.
We define an error benchmark
\[
\delta_\sigma := \sup \{ \| u - v \| : \dist(u, \cM), \dist(v,\cM) \le \sigma, P_W u = P_W v \} .
\]
This quantity can take in to account errors from model bias, and we have the following result.
\begin{theorem} \label{thm:rec_err}
Given a $\sigma$-admissible family of affine reduced spaces $(V_k)_{k=1}^K$, the estimator based on the surrogate selection \eqref{eq:surr_sel} has worst-case error bounded above by
\begin{equation} \label{eq:surr_err}
\max_{u \in \cM} \| u - u^*(P_W u) \| \le \delta_{\kappa \sigma},
\end{equation}
where $\kappa = R / r$ depends only on the uniform bounds of the surrogate distance.
\end{theorem}
The proof is detailed in Theorem 3.2 of \cite{cohen_nonlinear_2020}. Note that even given some optimal nonlinear reconstruction algorithm, our best possible error would be $\delta_0$, and it is not bounded from below by $d_n(\cM)$. We remark also that $\delta_\sigma \ge \sigma$. 

\section{Affine elliptic operators}
 Say the operator $A(y)$ in \eqref{eq:pde} is uniformly bounded in $y$ with uniformly bounded inverse, that is for some $0 < r \le R < \infty$ we have
\begin{equation} \label{eq:A_bounds}
\| A(y) \|_{\cH \to \cH^\prime} \le R \quad \text{and} \quad \| A(y)^{-1} \|_{\cH^\prime \to \cH} \le r^{-1} \quad y \in \cY.
\end{equation}
Then we can show that for any $v \in \cH$ the residual of the PDE,
\[
\cR(v, y) := \| A(y) v - f(y) \|_{\cH^\prime},
\]
satisfies the uniform bound $r \| v - u(y) \| \le \cR(v, y) \le R \| v - u(y) \|$. If we define the following,
\[
\cS(v, \cM) = \min_{y \in \cY} \cR(v, y),
\]
then we have arrived at a surrogate distance that satisfies the uniform bounds of \eqref{eq:surr_bounds}. Using this surrogate in the selection \eqref{eq:surr_sel} to define $u^*(w)$, we have a nonlinear reconstruction algorithm with the error guarantees of \eqref{eq:surr_err}.

We now make the further assumption that the operator $A(y)$ and source term $f(y)$ have affine dependence on $y$, that is 
\[
A(y) = A_0 + \sum_{j=1}^d A_j \quad \text{and} \quad f(y) = f_0 + \sum_{j=1}^d y_j f_j.
\]
The residual can be calculated using representers in $\cH$. We define $e_j$ as member of $\cH$ that satisfies
\begin{equation} \label{eq:riesz}
\< e_j, z \> = \< A_j(y) v - f_j , z \>_{\cH^\prime, \cH} \quad \text{for all}\quad z \in \cH,
\end{equation}
and now write $e(y) := e_0  -  \sum_{j=1}^d y_j e_j$ to denote the representer of the overall residual problem. 
The residual is equal to $\cR(v, y) = \| e(y) \|$, and determining the surrogate distance is a quadratic minimisation problem
\[
\cS(v, \cM)^2 = \min_{y \in \cY} \| e(y) \|^2 = \min_{y \in \cY} \big \| e_0 + \sum_{j=1}^d y_j e_j \big \|^2 .
\]
This leads to a constrained least squares problem that can be solved using standard optimisation routines, using the values $\< e_i, e_j \>$ for $0 \le i,j \le d$.

\subsection{Finite element residual evaluation}

In practice these calculations take place in a finite element space $\cH_h$ that is the span of polynomial elements on a triangulation $\cT_h$ of width $h>0$. In this setting the residual is
$\cR_h(v, y) = \| e_h(y) \| $
where $e_h(y) =  e_{0,h} -  \sum_{j=1}^d y_j e_{j,h}$
and the $e_{j,h} \in \cH_h$ satisfy the variational problem
\begin{equation} \label{eq:h_var}
\< e_{j,h}, z \> = \< A_j(y) v - f_j, z \>_{\cH^\prime, \cH} \quad \text{for all}\quad z \in \cH_h .
\end{equation}
Naturally we define $\cS_{h}(v, \cM) := \min_{y \in \cY} \cR_h(v,y)$. 
Note that when we subtract \eqref{eq:riesz} from \eqref{eq:h_var} we obtain
$\< e_h(y) - e(y) , z \> = 0$ for all $z \in \cH_{h}$, meaning that $e_h(y) - e(y) \perp \cH_h$.

We write $y^* = \argmin_{y \in \cY} \cR(v, y)$ to be the minimiser selected in $\cS(v, \cM)$, and $y^*_h$ the equivalent for $\cS_h(v, \cM)$. In general $y^* \neq y^*_h$, but we have
\begin{align}
| \cS(v, \cM) - \cS_h (v, \cM) | 
&= | \cR(v, y^*) - \cR_h(v, y^*_h) |  \nonumber \\
&\le | \cR(v, y^*) - \cR_h(v, y^*) | + | \cR(v, y^*_h) - \cR_h(v, y^*_h) | \nonumber \\
&\le \| e(y^*) - e_h(y^*) \| + \| e(y^*_h) - e_h(y^*_h) \|, \label{eq:S_err}
\end{align}
where in the last step we have used the fact that
\[
\left | \cR(v, y) - \cR_h(v, y) \right | 
= \big | \| e(y) \| - \| e_h(y) \| \big |
\le \left \| e(y) - e_h(y) \right\| .
\]
Thus the convergence of $\cS_h$ to $\cS$ depends on the finite element convergence of solutions $e_h$ to $e$. This convergence is determined by the regularity of $e(y)$, which depends on the smoothness of the data $A(y) v - h$ and the so called Riesz lift implied in the variational problem \eqref{eq:riesz}.


Recall that $A(y)$ in \eqref{eq:pde} is a second order symmetric elliptic operator. If we assume homogeneous Dirichlet boundary conditions on $\partial D$ and that $f(y) \in L^2(D)$, then a natural choice for our ambient space is the Sobolev space $\cH = H^1_0(D)$, with $\| \cdot \| = \| \cdot \|_{H^1_0(D)}$.
In this setting the solutions $e(y)$ of \eqref{eq:riesz} are the weak solutions of the Poisson problem with homogeneous Dirichlet boundary conditions, 
\begin{equation} \label{eq:poisson}
\nabla_x^2 e(y) = A(y) v - f(y) \text{ on } D \text{ with } e(y) = 0 \text{ on } \partial D, 
\end{equation}
where we have written $\nabla_x^2$ to denote the Laplacian in the spatial variables, recalling that we have the unwritten spatial dependence in $e(y) = e(x, y)$.

In our surrogate model selection, $e(y)$ depends on $u^*_k(w)$. This estimator $u^*_k(w)$ lies in $V_k \oplus W$. As $V_k$ is the span of some solutions selected from $\cM_k$, the smoothness of $u^*_k(w)$ will depend on the smoothness of all solutions of the PDE problem and of the measurement functionals $\omega_i$. 

The convergence of $e_h$ to $e$ in the weak form of \eqref{eq:poisson} is well known in a wide variety of settings, for example we have the classical result
\[
\| e_h(y) - e(y) \| \le c h \| A(y) v - f(y) \|_{L^2(D)},
\]
which would be applicable in a wide variety of situations. Under these circumstances we thus have that 
$| \cS(v, \cM) - \cS_h(v, \cM) | \sim h$.

\subsection{Coarse surrogate evaluation}

Say we construct a family of reduced spaces, $(V_k)_{k=1}^K$, where each family $V_k = \overline u_k + \spn(u^1_{h^\prime},\ldots,u^{n_k}_{h^\prime})$, and the solutions $u^i_{h^\prime}$ are numerically calculated with respect to a triangulation $\cT_{h^\prime}$ with mesh-width $h^\prime$. Our estimators $u^*_k(w)$ will  be in $\cH_{h^\prime}$.

We may consider $\cT_{h^\prime}$ to be our {\em fine} mesh, and nominate another triangulation $\cT_{h}$ with $h^\prime < h$ to be a {\em coarse} mesh. We can use this coarse mesh to compute the surrogate distance, noting that $\cS_h(v, \cM)$ may be orders of magnitude faster to compute than the fine mesh equivalent. This will necessarily introduce an inaccuracy in the surrogate selection, however we maintain the high fidelity of the fine mesh reduced space approximations $V_k$. 

If $\cT_{h^\prime}$ is a fine mesh that contains $\cT_h$ in the sense that $\cH_h \subset \cH_{h^\prime}$, 
then the variational problem \eqref{eq:h_var} straightforward as we can calculate the inner-product $\< A(y) v - f(y), z\>$ based on the known relationships between the basis elements of $\cH_h$ and $\cH_{h^\prime}$.
Furthermore the size of the error $\| e_h(y) - e_{h^\prime}(y) \|$, and hence $| \cS_h(v, \cM) - \cS_{h^\prime}(v, \cM) |$, is bounded above by a constant times $\| e(y) - e_h(y) \|$ through the same reasoning as in \eqref{eq:S_err}.

\section{Numerical tests} \label{sec:numerics}

We present two separate studies of the surrogate as evaluated in finite element spaces. On the unit square $D = [0,1]^2$ we consider the PDE
\[
- \nabla_x \cdot( a(x,y) \nabla_x u(x,y)) = 1 \quad \text{with} \quad u(x,y)=0 \text{ on } \partial D,
\]
where our parametric diffusivity field is given by $a(x,y) = 1 + \sum_{j=1}^{16} c_j y_j \Chi_{D_j}(x)$. Here the $D_j$ denote squares of side length $1/4$ that subdivide the unit square in to 16 portions,
and $\Chi_{D_j}$ is the indicator function on $D_j$. The parameter range is the unit hypercube $\cY = [-1, 1]^{16}$, and the coefficients are $c_j = 0.9 j^{-1}$ or $c_j = 0.99j^{-1}$, meaning that $a(x,y) > 0$, but when $c_j=0.99j^{-1}$ the $a(x,y)$ can become closer to 0 hence the PDE problem can lose ellipticity.

We perform space discretisation by the Galerkin method using $\mathbb{P}_1$ finite elements to produce solutions $u_{h^\prime}(y)$, with a fine triangulation on a regular grid of mesh size $h^\prime=2^{-7}$. In the tests that follow, we will evaluate $\cS_h(v, \cM)$ with coarse meshes of size $h=2^{-s}$ with $s = 2,\ldots, 6$.

We generate training sets to compute the reduced models, and test sets on which we test the reconstruction algorithms. 
The training set $\wt\cM^{\mathrm{tr}}$ is taken to be the collection of PDE solutions for $N_{\mathrm{tr}}=1000$ random samples 
 $\wt \cY^{\mathrm{tr}}=\{ y^{\mathrm{tr}}_j \}_{j=1,\ldots,N_\mathrm{tr}}$ drawn independently and uniformly on $\cY=[-1,1]^{16}$. The test set $\wt\cM^{\mathrm{te}}$ is created from $N_{\mathrm{te}}=100$ independent parameter samples that are distinct from the training set samples.


In our test the measurement space $W$ is given by a collection of $m=8$ measurement functionals $ \<\omega_i, u\> = \ell_i(u) = |B_i|^{-1} \int u \, \Chi_{B_i}$ that are local averages in a small area $B_i$ which are boxes of width $2h = 2^{-6}$, each placed randomly in the unit square.


In our first test we compute a single reduced space $V$ using the greedy procedure on $\wt \cM^{\mathrm{tr}}$. For each test candidate $u\in\wt\cM^{\mathrm{te}}$ we calculate $u^*(P_W u)$. In Figure \ref{fig:surr_err} we plot the difference between the coarse and fine surrogate distances $| \cS_h(u^*, \cM) -  \cS_{h^\prime}(u^*, \cM) |$. We note that this error is significantly dominated by the calculated value of the estimator error \eqref{eq:err_est}, $\sigma_{\mathrm{est}} := \mu(V, W) \, \e_{\mathrm{est}}(V, \cM)$
as we found $\sigma_{\mathrm{est}}\approx 10^0$ for $c_j = 0.9j^{-1}$ and  $\sigma_{\mathrm{est}}\approx 10^1$ for $c_j = 0.99j^{-1}$. This dominates the errors presented in Figure \ref{fig:surr_err} significantly. 

Furthermore, we see a linear relationship in the right hand in Figure \ref{fig:surr_err}. The linear best fit of $\log (\mathrm{avg}_{\wt\cM^{\mathrm{te}}} | \cS_h - \cS_{h^\prime} |)$ and $\log(h) = -s$ has slope $1.51$, which implies that on average 
\[
| \cS_h - \cS_{h^\prime} | \sim h^{1.51}
\]

\begin{figure}
\centering   
  \includegraphics[width=\linewidth]{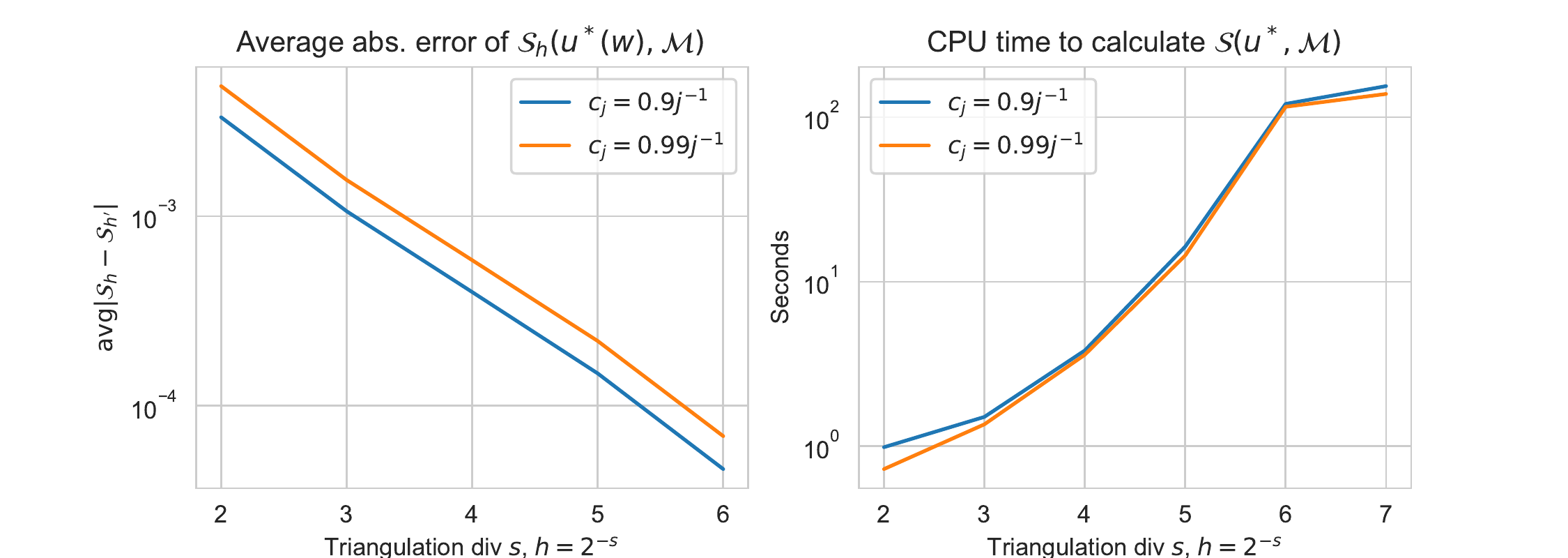}
  \vspace{-20pt}
  \caption{\small Left, average absolute error of the surrogate $\cS_h$, right, the CPU wall time for the computation of $\cS_h$ for all $N_{\mathrm{te}}$ test points}
  \vspace{-4pt}
  \label{fig:surr_err}
\end{figure}


For the second test we examine the impact of using the coarse surrogate $\cS_h(u^*_k, \cM)$ for model selection. We build the $\sigma$-admissible families as outlined in \S\ref{sec:nonlinear} using the greedy splitting of $\cY = [-1,1]^{16}$. Note that at each point the cells $\cY_k$ are rectangular, which we can split in half in coordinate directions.
We split the parameter space 7 times, resulting in $K=8$ local reduced spaces.

For each test candidate $u \in \wt \cM^{\mathrm{te}}$ we have $K$ possible reconstructions $u_1^*(w),\dots u_K^*(w)$. We use the coarse surrogate in the model selection \eqref{eq:surr_sel}, writing $k^*_h(w) = k^*_h(P_W u)$ to make the dependence on $h$ and $w$ clear, and we inspect how often it agrees with the {\em fine} selection $k^*_{h^\prime}(w)$ for all test points $u$. We can also compare compare this to and the ``true'' selection $k^*_{\mathrm{true}}(u)$ for which $u \in \cM_{k^*_{\mathrm{true}}}$.

Table \ref{tab:k_h} demonstrates that $k^*_h(w)$ agrees with $k^*_{h^\prime}$ the vast majority of the time from $h=2^{-4}$ onwards, in both cases for $c_j$. We also see that the fine selection $k^*_{h^\prime}$ agrees with the true selection 77 times out of 100 for $c_j = 0.9j^{-1}$ and 64 times for $c_j = 0.99j^{-1}$, that is, it picks the estimator $u^*_k(w)$ that is trained on the portion of manifold $\cM_k$ that $u$ originated from. 
Out of interest we also plot the histogram of selections $k^*_{h^\prime}$ and $k^*_{\mathrm{true}}(u)$, recalling that they are a number in $1,\ldots, 8$. We see broadly similar patterns in the reduced model selection. 
Given the CPU time savings that we see in Figure \ref{fig:surr_err}, we conclude that model selection through a coarse surrogate distance is a worthwhile numerical strategy.

\begin{table}
\centering\footnotesize

\begin{tabular}{| c | c | c | c | c | c | c | c | c | c | c |}
 \multicolumn{1}{c}{} & \multicolumn{5}{ | c | }{$c_j=0.9 j^{-1}$}  & \multicolumn{5}{ |c | }{$c_j=0.9 j^{-1}$}  \\
 \hline
Mesh width $h$ & $2^{-3}$ & $2^{-4}$ & $2^{-5}$ & $2^{-6}$ & $2^{-7}$ & $2^{-3}$ & $2^{-4}$ & $2^{-5}$ & $2^{-6}$ & $2^{-7}$ \\
\hline
$\# \{ k^*_h(w) = k^*_{h^\prime}(w) \}$ & 97 & 100 & 100 & 100 & 100  & 88 & 94 & 96 & 98 & 100  \\
$\# \{  k^*_h(w) = k^*_{\mathrm{true}}(u) \}$ & 74 & 77 & 77 & 77 & 77 & 58 & 59 & 61 & 65 & 64 \\
\hline
\end{tabular}
\caption{\small Agreement of coarse model selection out of $N^{\mathrm{te}} = 100$ test points $u \in \wt\cM^{\mathrm{te}}$. \label{tab:k_h}}
\end{table}

\begin{figure}
\centering   
  \includegraphics[width=\linewidth]{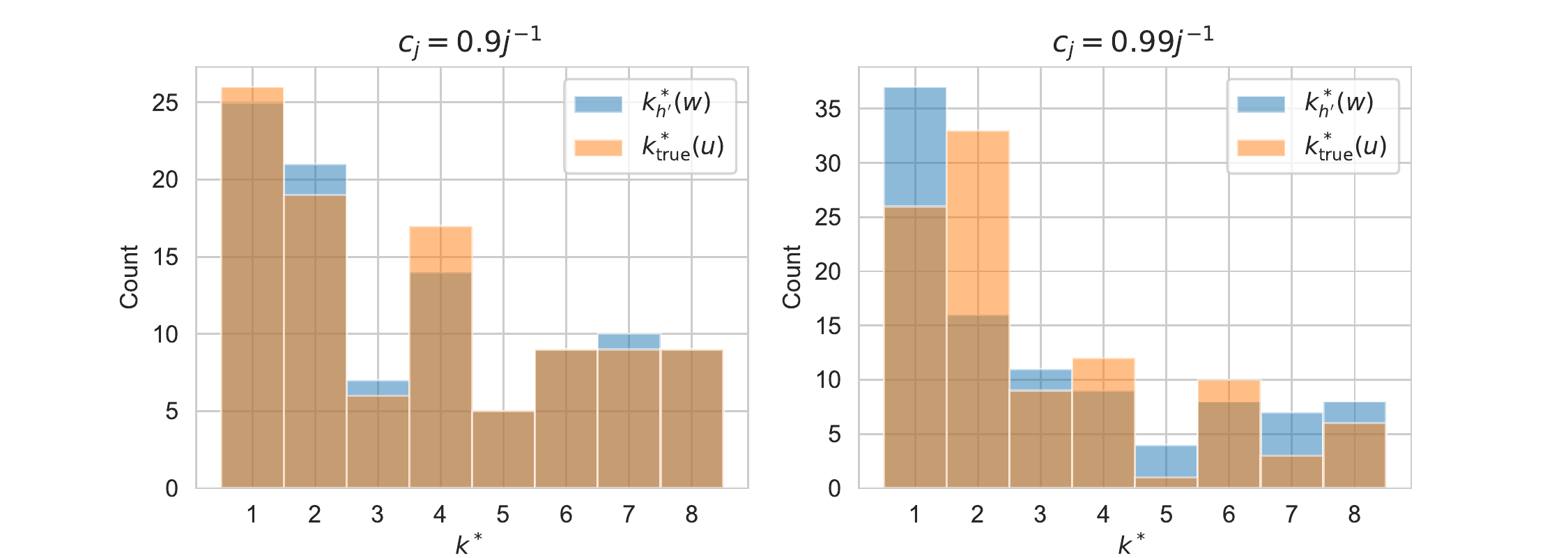}
  \vspace{-20pt}
  \caption{\small A histogram of the number of fine surrogate selections $k^*_{h^\prime}$ and true reduced model selection $k^*_{\mathrm{true}}(u)$.}
  \vspace{-4pt}
  \label{fig:surr_err}
\end{figure}

{
\bibliography{multires_refs}}

\begin{thebibliography}{1}

\bibitem{binev_convergence_2011}
Peter Binev, Albert Cohen, Wolfgang Dahmen, Ronald DeVore, Guergana Petrova,
  and Przemyslaw Wojtaszczyk.
\newblock Convergence {Rates} for {Greedy} {Algorithms} in {Reduced} {Basis}
  {Methods}.
\newblock {\em SIAM Journal on Mathematical Analysis}, 43(3):1457--1472,
  January 2011.
\newblock Publisher: Society for Industrial and Applied Mathematics.

\bibitem{binev_data_2017}
Peter Binev, Albert Cohen, Wolfgang Dahmen, Ronald DeVore, Guergana Petrova,
  and Przemyslaw Wojtaszczyk.
\newblock Data assimilation in reduced modeling.
\newblock {\em SIAM/ASA Journal on Uncertainty Quantification}, 5(1):1--29,
  2017.

\bibitem{cohen_nonlinear_2020}
Albert Cohen, Wolfgang Dahmen, Olga Mula, and James Nichols.
\newblock Nonlinear reduced models for state and parameter estimation.
\newblock {\em arXiv:2009.02687 [cs, math]}, November 2020.
\newblock arXiv: 2009.02687.

\bibitem{cohen_reduced_2020}
Albert Cohen, Dahmen Wolfgang, Ronald DeVore, and James Nichols.
\newblock Reduced basis greedy selection using random training sets.
\newblock {\em ESAIM: Mathematical Modelling and Numerical Analysis}, January
  2020.

\bibitem{hesthaven_certified_2016}
Jan~S. Hesthaven, Gianluigi Rozza, and Benjamin Stamm.
\newblock {\em Certified {Reduced} {Basis} {Methods} for {Parametrized}
  {Partial} {Differential} {Equations}}.
\newblock {SpringerBriefs} in {Mathematics}. Springer International Publishing,
  2016.

\bibitem{maday_parameterized-background_2015}
Yvon Maday, Anthony~T. Patera, James~D. Penn, and Masayuki Yano.
\newblock A parameterized-background data-weak approach to variational data
  assimilation: formulation, analysis, and application to acoustics.
\newblock {\em International Journal for Numerical Methods in Engineering},
  102(5):933--965, May 2015.

\end{thebibliography}
\bibliographystyle{plain}

\end{document}